\documentclass[10pt]{siamltex}

\usepackage{amsfonts}
\usepackage{amssymb}
\usepackage{times}
\usepackage{graphicx}

\usepackage{color}
\def\Image{\hbox{\rm Image}}
\def\R{\mathbb R}

\def\intx{\int\limits_{x_{i-1}}^{x_i}}
\def\inty{\int\limits_{\tilde x_{i}}^{\tilde x_{i+1}}}

\def\R{\mathbb{R}}

\def\Image{\mathop\mathrm{Image}}

\def\<{\langle}
\def\>{\rangle}

\title{A simple, fast and stabilized flowing finite volume method
for solving general curve evolution equations}

\author{Karol Mikula
\thanks{Department of Mathematics, Slovak University of
Technology, Rad\-lin\-sk\'eho 11, 813 68 Bratislava, Slovak Republic
({\tt mikula@math.sk}).}
\and
{Daniel \v Sev\v covi\v c}
\thanks{Department of Applied Mathematics and Statistics,
Faculty of Mathematics, Physics \& Informatics, Comenius University, 842 48
Bratislava, Slovak Republic ({\tt sevcovic@fmph.uniba.sk}).}
\and
{Martin Bala\v zovjech}
\thanks{Department of Mathematics, Slovak University of
Technology, Rad\-lin\-sk\'eho 11, 813 68 Bratislava, Slovak Republic
({\tt balazoviech@math.sk}).
\newline \newline This work was supported by grants: 
VEGA 1/0269/09, APVV-0351-07, APVV-RPEU-0004-07 
(K.Mikula and M.Bala\v zovjech) and APVV-0247-06 (D.\v{S}ev\v{c}ovi\v{c}).
}
}

\begin{document}

\maketitle

\begin{abstract}A new simple Lagrangian method with favorable stability and
efficiency properties  for computing general plane curve evolutions is presented.
The method is based on the  flowing finite volume discretization of the 
intrinsic partial differential  equation for updating the position vector of evolving family of plane curves. 
A curve can be evolved in the normal direction by a combination of fourth order terms related to the intrinsic Laplacian of the curvature,  second order terms related to the curvature, first order terms related to
anisotropy and by a given external velocity field. The evolution is numerically stabilized by an asymptotically uniform tangential redistribution of grid points yielding the first order intrinsic advective terms in the governing system of equations. By using a semi-implicit in time discretization it can be numerically approximated by a solution to linear penta-diagonal systems of equations (in presence of the fourth order terms) or tri-diagonal systems (in the case of the second order terms). Various numerical experiments of plane curve evolutions, including, in particular, nonlinear, anisotropic and regularized backward curvature flows, surface diffusion and Willmore flows, are presented and discussed.
\end{abstract}



\pagestyle{myheadings}
\thispagestyle{plain}

\section{Introduction}

The main purpose of this paper is to propose a simple, fast and stable 
Lagrangian method for computing evolution of 
closed smooth embedded plane curves  driven by a normal velocity of the form
$\beta(\partial^2_s k, k, \nu, x)$ which may depend
on the intrinsic Laplacian  $\partial^2_s k$ of the curvature 
$k$, on the curvature $k$ itself, on the tangential angle $\nu$ and the curve position vector $x$. 
We shall restrict our attention to the following form of the normal velocity:
\begin{equation}
\beta = -\ \delta\ \partial^2_s k  + b(k,\nu) + F(x).
\label{geomrov}
\end{equation}
Here $\delta\ge 0$ is a constant and $b=b(k,\nu)$ is a smooth function 
satisfying $b(0,.)=0$.  If $\delta>0$ then there is no constraint on the monotonicity of $b$ with respect to $k$. On the other hand, if $\delta=0$, then we shall assume the function $b$ is strictly increasing with respect to the curvature.

There are many interesting evolutionary models in various applied fields of science, technology and engineering that can be described by geometric equation (\ref{geomrov}). For example, putting $\delta=0$ we obtain the normal velocity $\beta=b(k,\nu) + F(x)$ representing the well-known anisotropic mean curvature flow arising in the motion of material interfaces during solidification and in affine invariant shape analysis (see e.g. \cite{AG1,GH,Gr,D2,D3,ST,MK,M}). The term $F(x)$ represents an external driving force, like e.g., a given velocity field projected to the normal vector of a curve, or any other constant or scalar function depending on the current curve position. It drives the curve in the inner (if $F(x)>0$) or outer (if $F(x)<0$) normal direction. 
In the case $\delta=1$ two well-known examples  arise when studying a motion of the so-called elastic curves. It is 
the surface diffusion, in the case $b=0$, and the Willmore flow  
where $b = -\frac12 k^3$. The surface diffusion is often used in computational fluid dynamics and material sciences, where the encompassing area of interface should be preserved. 
The case when $b = -\frac12 k^3$ arises from the model of the Euler-Bernoulli elastic rod -- 
an important problem in structural mechanics \cite{CT,E,DKS}. 
The evolutionary models having the normal velocity of the form 
(\ref{geomrov}) are often adopted in image segmentation where elastic and 
geodesic curves are used in order to find image objects in an  automatic 
way \cite{KWT,CKS2,KKOTY2,MS_CVS,MS_AA}. 
By our method we are able to handle a regularized backward mean curvature 
flow in which $b$ is a decreasing function of the curvature $k$ 
like e.g. $b(k,\nu)=-k$. We regularize the backward mean curvature 
flow by adding a small fourth order 
diffusion term  $0<\delta\ll 1$. 
To our knowledge, first experiments of this kind are presented in this paper. 

The main idea of our approach is based on accompanying the geometric equation
(\ref{geomrov}) by a stabilizing tangential velocity and in rewriting it into a form of intrinsic an
partial differential equation (PDE) for the curve position vector. 
The resulting PDE contains fourth, second and first order spatial differential terms that 
are approximated by means of the flowing finite volume method \cite{MS2}. For time
discretization we follow semi-implicit approach leading  to a solution to a 
linear system of equations at each time level. That can be done efficiently,
and, due to tangential stabilization, 
we hope that this direct Lagrangian method can be considered 
as an efficient counterpart to the well-known level-set based methods 
for description of the curve evolution discussed in  
\cite{OS,Se2,Osherbook,HMSg1,FM1,FM2}.

The paper is organized as follows: in the next section we derive the intrinsic PDE for description of a family of plane curves. Then we present our numerical approximation scheme. Finally, we discuss various numerical experiments showing applicability of our approach. 

\section{Governing equations}
An immersed regular plane  curve $\Gamma$ can be parameterized by a smooth function 
$x:S^1\to \R^2$, i.e. $\Gamma=\{ x(u), u\in S^1 \}$ having a strictly positive local length term 
$g=|\partial_u x| >0$. Taking into account periodic boundary conditions at 
$u=0,1$ we can identify the circle $S^1$ with the interval $[0,1]$. The unit arc-length parameterization is denoted by 
$s$. Clearly, $\mbox{d} s = g\, \mbox{d} u$. For the tangent vector we have $\vec T = \partial_s x$
and we can choose the unit inward normal vector $\vec N$  such that 
$\mbox{det}(\vec T,\vec N) =1$.
The tangential angle $\nu=\arg(\vec T)$,
$\vec T = (\cos \nu, \sin \nu )^\top$ and $\vec N = (-\sin \nu, \cos \nu )^\top$.  

Let a regular smooth initial curve $\Gamma^0=\Image(x^0)=\{x^0(u), u\in[0,1]\}$ be given.
We shall represent a family of plane curves 
$\Gamma^t=\Image\ (x(.,t))=\{x(u,t), u\in[0,1]\}$ that evolves according to the geometric equation  (\ref{geomrov}) by its position vector $x$ satisfying the following equation:
\begin{equation}
\partial_t x = \beta \vec N + \alpha \vec T\,.
\label{geomrov2}
\end{equation} 
It is well-known that the presence of any tangential velocity functional $\alpha$ in (\ref{geomrov2}) has no impact on the shape of
evolving curves. However, it may help to redistribute points along the curve. As a consequence, it
can  significantly stabilize numerical computations. The reader is referred to papers
\cite{Hou1,K2,MS1,MS2,MS3,MS_CVS,MS_ALG,De,BGN1,BGN2,BMOS,Y,SY} for detailed discussion on how a suitable tangential stabilization can prevent a numerical solution from forming various undesired singularities. We will specify our choice of a tangential velocity $\alpha$ later. Since
\begin{eqnarray}
\partial_s^4 x && = \partial_s^3 \vec T = \partial_s^2 (k \vec N) = 
\partial_s^2 k \vec N + 2 \partial_s k \partial_s \vec N + k\partial_s^2
\vec N =
\partial_s^2 k \vec N - 2 (\partial_s k) k \vec T - k\partial_s(k\vec T)
\nonumber\\
&&=\partial_s^2 k \vec N -3 k (\partial_s k) \vec T - k^2 \partial_s \vec T =
\partial_s^2 k \vec N -\frac {3}{2}\partial_s (k^2) \partial_s x - k^2
\partial_s^2 x,
\nonumber
\end{eqnarray} 
we have 
\begin{equation}
(-\partial_s^2 k) \vec N = -\partial_s^4 x -k^2
\partial_s^2 x - \frac {3}{2}\partial_s (k^2) \partial_s x.
\end{equation}
Let us define the following auxiliary functions:
\begin{equation}
\phi(k,\nu)=-\delta k^2 + c(k,\nu),\ \ \  c(k,\nu)= b(k,\nu)/k
\end{equation}
and 
\begin{equation}
v(k,\nu)=\frac {3}{2}\delta\ \partial_s (k^2)+\partial_s
\phi(k,\nu)-\alpha \,.
\label{rovnicaV}
\end{equation}
Since the function $b$ is assumed to be smooth and $b(0,\nu)=0$ the function 
$c(k,\nu)=b(k,\nu)/k$ is smooth as well. Using the Frenet formulae we obtain
$b(k,\nu) \vec N = c(k,\nu) \partial_s^2 x$ and 
$\phi(k,\nu)\partial_s^2 x = \partial_s(\phi(k,\nu)\partial_s x)-
\partial_s\phi(k,\nu)\ \partial_s x$. Hence, for the normal velocity $\beta$ of the form 
(\ref{geomrov}) we end up with the following higher-order intrinsic PDE for the position vector $x=x(s,t)$:
\begin{equation}
\partial_t x + v(k,\nu)\partial_s x = \delta(-\partial_s^4 x) +
\partial_s(\phi(k,\nu)\partial_s x) + F(x) 
\vec N(\nu).
\label{geomrov3}
\end{equation}
Numerical approximation of a solution $x$ to the above PDE forms the basis of our direct Lagrangian approach.

It is known (see e.g.  \cite{MS2})  that a family of plane curves $\Gamma^t = \Image(x(., t)), t\in [0,T)$, that evolves according to (\ref{geomrov2}) can be also represented  by a solution 
to the following system of intrinsic parabolic-ordinary differential
equations:
\begin{eqnarray}
&&\partial_t k =\partial^2_s \beta +\alpha \partial_s k +k^2\beta,
\label{rovnice1} \\
&&\partial_t \nu = \partial_s\beta
     + \alpha k,
\label{rovnice2} \\
&&\partial_t g  = -g k\beta  + g\partial_s\alpha .\label{rovnice3}
\end{eqnarray}
In \cite{MS3,MS_CVS,MS_ALG} the system (\ref{rovnice1})-(\ref{rovnice3}) was solved numerically for the curvature $k$, tangent angle $\nu$ and the local length $g$. Knowing these quantities one can 
reconstruct the curve evolution. 
Asymptotically uniform redistributions for the second order flows 
with driving force (see \cite{MS3}) and for the fourth order flows (see 
\cite{MS_ALG}) were also proposed. 

In this paper we follow a different approach when compared to \cite{MS3,MS_ALG}. It is much faster and  simpler 
from computational point of view. We do not solve the system (\ref{rovnice1})-(\ref{rovnice3}), but
the position vector equation (\ref{geomrov3}) is directly  numerically discretized. 
On the other hand, from the analytical point of view, the system (\ref{rovnice1})-(\ref{rovnice3}) 
describes evolution of useful geometric quantities that can be utilized in
designing proper tangential velocities for stabilization of numerical computations. 
For example, in the case of convex curves, it is sufficient to
solve just equation (\ref{rovnice1}) on the fixed parameter interval given by
the range of $\nu$. Then no grid point redistribution is necessary
\cite{MK,M}. Equation (\ref{rovnice2}) was used recently in designing
the tangential velocities corresponding to the so-called crystalline
curvature flow, where $\partial_t \nu=0$ (c.f. \cite{Y,SY}). 
In \cite{MS3} the authors proposed and analyzed a new type of tangential redistribution referred to as the 
{\it asymptotically uniform grid points redistribution}. Here we follow a similar idea but adopt it in a rather different way. More precisely, in order to compute the corresponding tangential velocity $\alpha$ (depending on the curvature,tangential angle and local length) we use information from a solution to (\ref{geomrov3}) only, i.e. we do not solve the system (\ref{rovnice1})-(\ref{rovnice3}).

We will see that such a simple and straightforward 
approach yields a stable, fast and precise 
solution to our problem. 
The corresponding numerical scheme is efficiently stabilized by an 
appropriate choice of the tangential velocity 
term $\alpha$ in (\ref{rovnicaV}).
Let us note that other known tangential velocities, 
like the one preserving relative local length \cite{Hou1,K2,MS2}, 
locally diffusive redistribution
\cite{De,MS_CVS,MS_ALG}, crystalline curvature redistribution \cite{Y} 
or curvature adjusted redistribution
\cite{SY,PB} can be also incorporated straightforwardly 
to the numerical scheme
by corresponding change of the term $\alpha$ 
entering equation (\ref{rovnicaV}).

Let us briefly outline the basic idea behind the asymptotically uniform grid points redistribution.
Let us denote by $L_t=\int_{\Gamma^t} ds=\int_0^1 g(u,t)\,du$ 
the length of an evolving curve $\Gamma^t$. Integrating equation
(\ref{rovnice3}) along the curve and taking into account 
periodicity of $\alpha$ at $u=0,1$ we obtain
\begin{equation}
\frac{d L}{dt}  + \langle k \beta \rangle_{\Gamma} L =0,
\label{length-eq}
\end{equation}
where $\langle k \beta \rangle_{\Gamma} =  \frac{1}{L}\int_\Gamma k \beta \, ds$ 
denotes the curve average of the quantity $k \beta$ over a curve $\Gamma$.

When designing the tangential velocity $\alpha$, it is worth to study the time evolution of
the quantity $g/L$. 
From a numerical point of view it corresponds to the 
local grid point distances divided by the averaged local distance $L/n$, where
$n$ is the number of grid points. 
We refer the reader to content of the next section for details. 
In that sense, it represents a deviation of local grid point distances 
from being uniformly distributed. Moreover, if we define the quantity $\theta = \ln (g/L)$ and
we take into account equations (\ref{rovnice3}) and 
(\ref{length-eq}) we conclude
\begin{equation}
\partial_t \theta + k\beta -\langle k\beta \rangle_\Gamma 
= \partial_s \alpha\,.
\label{thetarovnica}
\end{equation}
From (\ref{thetarovnica}) we can observe that, by an appropriate choice of 
$\partial_s \alpha$, we can control behavior of $\theta$ and, subsequently, of the ratio $g/L$ also. 
Our choice of $\alpha$ is based on the following particular setup (see \cite{MS3})
\begin{equation}
\partial_s \alpha = k\beta -\langle k\beta \rangle_\Gamma
+ \left(  e^{-\theta} -1\right) \omega(t)
\label{uniformredis}
\end{equation}
where $\omega\in L^1_{loc}([0,T_{max}))$ and $T_{max}$ is the maximal
existence time of evolving curve. 

The most simple choice $\omega(t)\equiv 0$ 
yields $\partial_t\theta=0$. Hence
\[
\frac{g(u,t)}{L_t} = \frac{g(u,0)}{L_0} \quad \hbox{for any} \ u\in S^1,
\ t\in [0,T_{max}),
\]
and we obtain the so-called tangential redistribution preserving the relative local length 
(cf. \cite{Hou1,K2,MS2}). On the other hand, assuming
\begin{equation}
\int_0^{T_{max}} \omega(\tau) \, d\tau = + \infty
\label{omegaintegral}
\end{equation}
then, by inserting  (\ref{uniformredis}) into (\ref{thetarovnica})  
and solving the corresponding ordinary differential equation 
$\partial_t\theta =  \left(  e^{-\theta} 
-1\right) \omega(t),$ 
we obtain $\theta(u,t) \to 0$ as $t\to T_{max}$ and hence
\[
\frac{g(u,t)}{L_t} \to 1  \quad \hbox{as}\ t\to T_{max} 
\quad \hbox{uniformly w.r. to} \ u\in [0,1] .
\]
It means that redistribution of grid points along a curve becomes {\it asymptotically  uniform} as $t$  approaches the maximal time of existence $T_{max}$. In the case when the family $\{\Gamma^t, t\in[0,T)\}$ shrinks to a point as 
$t\to T_{max}$, in order to fulfill (\ref{omegaintegral}), one can choose
$\omega(t) =  \kappa_2 \langle k\beta \rangle_{\Gamma^t}$ 
where $\kappa_2>0$ is a positive constant. By (\ref{length-eq}) we have
$\int_0^t \omega(\tau)\, d\tau 
= - \kappa_2 \int_0^t \ln L^\tau d\tau  = \kappa_2(\ln L^0 - \ln L_t) \to 
+\infty \quad \hbox{as}\ t\to T_{max}$ because $\lim_{t\to T_{max}} L_t = 0$. 
On the other hand, if the length $L_t$  is always away  from zero and $T_{max}=+\infty$ 
one can take $\omega(t) = \kappa_1$, where $\kappa_1>0$ is a positive constant. 
Summarizing, a suitable choice of the tangential velocity functional $\alpha$ that helps to redistribute grid points
uniformly along the evolved curve (in any case of shrinking, expanding or reaching an equilibrium curve shape) is given by a solution to the equation
\begin{equation}
\partial_s \alpha = k\beta -\langle k\beta \rangle_\Gamma
+ \left(L/g -1\right) \omega,\ \ \omega = \kappa_1 
+\kappa_2\langle k\beta\rangle_\Gamma, \ \ \alpha(0,t)=0
\,,
\label{alpha-nonlocal}
\end{equation}
where $\kappa_1,\kappa_2\ge 0, \kappa_1+\kappa_2>0,$ are given constants. 
In all computations to follow in Section 4, we
choose $\kappa_2 = 0$ which is sufficient when the computations end up
before an eventual extinction of a curve.
By the boundary condition imposed on $\alpha(0,t)$ we prescribe the
motion of the point $x(0,t)$ in the normal direction only. Therefore, the tangential velocity term $\alpha$ is determined uniquely.

It is worth to note that the speed of relaxation of 
the quantity $\theta=\ln(g/L)$ is controlled by the constant 
$\kappa_1>0$. We typically take $\kappa_1=O(1)$ in all our computations. More precisely, 
in experiments presented in section 4 we take $\kappa_1\in[3,10]$.   
This is due to the fact that the speed of the tangential velocity 
$\alpha$ defined as in (\ref{alpha-nonlocal}) is increasing with 
respect to the parameter $\kappa_1$. Therefore choosing considerably larger values 
of $\kappa_1$ would make the governing equation (\ref{geomrov3}) 
strongly convective dominant which may yield a necessity of small time 
steps in our numerical scheme. 

\section{Numerical approximation scheme}

A family of plane curves that evolves according to (\ref{geomrov}) is numerically represented by a family of "flowing" discrete plane points $x_i^j$ where the index $i=1,...,n,$ denotes space 
discretization and the index $j=0,...,m,$ stands for a discrete time
stepping. Assuming a uniform division of the time interval $[0,T]$ with a time step
$\tau=\frac{T}{m}$ and a uniform division of the fixed
parameterization interval $[0,1]$ with a step $h=1/n$, a discrete point
$x_i^j$ corresponds to $x(ih,j\tau)$. 
Due to periodic boundary conditions and smoothness of the evolved curve 
we have used the additional points
defined by  $x_{-1}^j=x_{n-1}^j$, $x_0^j=x_n^j$, 
$x_{n+1}^j=x_1^j$, $x_{n+2}^j=x_2^j$. 
The tangential velocity of a flowing node $x_i^j$ is denoted by $\alpha_i^j$.

The system of difference
equations corresponding to (\ref{alpha-nonlocal}) and (\ref{geomrov3})
will be constructed at each discrete time step $j\tau$ 
by using the flowing finite volume method proposed in \cite{MS2}.
First we solve (\ref{alpha-nonlocal}) for the tangential velocities
$\alpha_i^j$ and then  equation (\ref{geomrov3}) for the position vectors 
$x_i^j, i=-1, ..., n+2$. Remaining quantities
involved in (\ref{alpha-nonlocal}) and (\ref{geomrov3}), as the curvature,
tangential angle, local and total length of the curve 
are computed from the curve position vector $x^{j-1}$ from 
the previous time step $j-1$. 

In order to build our numerical scheme we construct the so-called 
{\it flowing finite volume} $[x_{i-1}^j,x_{i}^j]$ and  
also corresponding flowing dual 
volumes $[\tilde x_{i}^j,\tilde
x_{i+1}^j]$ where $\tilde x_i^j= (x_{i-1}^j+x_{i}^j)/2$. 
Then the approximate local lengths of flowing finite volumes 
$r_i^j$,  curvatures $k_i^j$ and tangential angles $\nu_i^j$ are given by
piecewise constant values in the flowing
finite volumes. Similarly, the other quantities  
$\alpha_i^j$, $x_i^j$ and approximate lengths of dual volumes 
$q_i^j$ are considered as piecewise constant
in flowing dual volumes. 

Let a discrete representation  of the evolving curve be given
at time level $j-1$ by discrete points $x_i^{j-1}, i=-1, ..., n+2$. 
At every new time level $j$ we compute
\begin{eqnarray}
r_i^j &=& |R_i|,\ R_i = (R_{i_1},R_{i_2})=x_i^{j-1}-x_{i-1}^{j-1}, 
\ \ i=0,...,n+2,\ q_i^{j}=\frac12 \left(r_{i}^{j}+r_{i+1}^{j}\right)\,,
\nonumber\\
k_i^j &=& {1\over 2 r_i^j}\hbox{sgn}(\det(R_{i-1}, R_{i+1}))
\arccos\left( {R_{i+1} . R_{i-1}\over r_{i+1}^j r_{i-1}^j}\right)\,,
i=1,...,n,\ k^j_0=k^j_n, \ k^j_{n+1}=k^j_1\,,
\nonumber\\
\nonumber\\
\nu_0^j &=& \arccos(R_{0_1}/r_0^j),\ \hbox{if}\ R_{0_2}\ge 0\,,\ \ \ 
\nu_0^j = 2\pi - \arccos(R_{0_1}/r_0^j), \ \hbox{if}\ R_{0_2}< 0\,,
\nonumber\\
\nu_i^j &=& \nu_{i-1}^j+r_i^j k_i^j,\ \ i=1,...,n,\ \ 
\nu_{n+1}^j=\nu_1^j+2\pi\,,
\nonumber\\
\nonumber\\
\beta_i^{j} &=& \frac{\delta}{r_i^{j}}
\left(
\frac{k_{i}^{j}-k_{i-1}^{j}}{q_{i-1}^{j}}-
\frac{k_{i+1}^{j}-k_i^{j}}{q_i^{j}}\right)+
b(k_i^{j},\nu_i^{j})+\frac{F(x_i^{j-1})+F(x_{i-1}^{j-1})}{2},
i=1,...,n,
\nonumber\\
\nonumber\\
L^{j} &=& \sum\limits_{l=1}^n r_l^{j},\ \ \ 
B^{j} = \frac{1}{L^{j}}
\sum\limits_{l=1}^n  r_l^{j}k_l^{j} \beta_l^{j}\,.
\nonumber
\end{eqnarray}

In order to compute the tangential velocity $\alpha$ we
integrate (\ref{alpha-nonlocal}) 
over the time varying flowing finite volume $\left[x_{i-1},x_i\right]$,
\[
\intx\partial_s\alpha \,\hbox{d}s = \intx
k\beta-\langle k\beta \rangle_\Gamma + \left(L/g - 1\right) \omega
\,\hbox{d}s\,.
\]
 
Hereafter we use the notation $\intx\psi \,\hbox{d}s$ for integral of the quantity $\psi$ over the curve arc $\widehat{x_{i-1},x_i}$.
Hence at any time level $t$ we have the relation 
for approximation of the difference $\alpha_i-\alpha_{i-1}$
\[
\alpha_i -\alpha_{i-1} \approx  r_i (k_i 
\beta_i- \langle k\beta \rangle_\Gamma) +\left(h L - r_i\right)\omega\,.
\]
Taking into account discrete time stepping in the previous relation we
obtain the following expression for up-dated values of the tangential velocity at
the $j$-th time level:
\begin{equation}
\alpha_i^j = \alpha_{i-1}^j + 
r_i^{j}(k_i^{j}\beta_i^{j} -
B^{j})+(h L^{j} - r_i^{j})\omega,\ i=1,...,n,\ \ \alpha_0^j=0.
\label{alpha-update}
\end{equation}
After computing tangential velocities we update diffusion and advection
terms using
\[
\phi_i^j = -\delta (k_i^j)^2 + 
c(k_i^j,\nu_i^j),\ \ i=1,...,n+1,
\]
\[
v_i^j=\frac{3\delta}{2}\frac{(k_{i+1}^j)^2-(k_i^j)^2}{q_i^j}+
\frac{\phi_{i+1}^j-\phi_{i}^j}{q_i^j}-\alpha_i^j,\ \ i=1,...,n\,.   
\]
In order to compute the curve position vector $x^j$ at a new time level we
integrate (\ref{geomrov3}) over the time varying flowing
dual volume
$\left[\tilde x_{i},\tilde x_{i+1}\right]$,
\[
\inty\partial_t x +v(k,\nu) \partial_s x\,\hbox{d}s = \inty 
-\delta\partial_s^4 x + \partial_s(\phi(k,\nu)\partial_s x)
+F(x)\vec N(\nu)
\,\hbox{d}s\,.
\]
Notice that the advective term $v(k,\nu)$ contains the tangential velocity $\alpha$ and 
partial derivatives of functions of $k$ and $\nu$. They  can be approximated
by constant difference quotients in flowing dual volumes. Therefore we can write
\begin{equation}
\ \ \ \ q_i\frac{\hbox{d} x_i}{\hbox{d}t} +v_i (\tilde x_{i+1}-\tilde x_i)=
\left[-\delta\partial^3_s x + \phi(k,\nu)\partial_s x 
\right]_{\tilde x_i}^{\tilde x_{i+1}}
+ q_i^j F(x_i)\vec
N\left(( \nu_{i}+\nu_{i+1})/2 \right).
\label{ode}
\end{equation}
Denote by $\partial^k_s x_i$ and $\partial^k_s \tilde x_i$ the approximation of the $k$-th arc-length derivative $\partial^k_s x$ at points $x_i$ and $\tilde x_i$, respectively. 
The differences at points $\tilde x_i$, $\tilde x_{i+1}$ 
in the above boundary integral terms are then naturally approximated by 
\[
\phi(k,\nu)\partial_s \tilde x_{i+1}- 
\phi(k,\nu)\partial_s \tilde x_i \approx
\phi_{i+1}^j\frac {x_{i+1}^j-x_{i}^j}{r_{i+1}^j}
- \phi_{i}^j\frac {x_{i}^j-x_{i-1}^j}{r_{i}^j}\,.
\]
The third order terms appearing in the boundary integral in (\ref{ode}) can be approximated as follows:
\begin{eqnarray*}
&&\partial_s^3 \tilde x_{i+1}-\partial_s^3
\tilde x_{i}\approx
\frac{\partial_s^2 x_{i+1}-\partial_s^2 x_{i}}{r_{i+1}}
-\frac{\partial_s^2 x_{i}-\partial_s^2 x_{i-1}}{r_{i}}
\\
&&\approx \frac{1}{r_{i+1}}\left(
\frac{\partial_s \tilde x_{i+2}-\partial_s \tilde x_{i+1}}{q_{i+1}}-
\frac{\partial_s \tilde x_{i+1}-\partial_s \tilde x_{i}}{q_{i}}
\right)
-
\frac{1}{r_{i}}\left(
\frac{\partial_s \tilde x_{i+1}-\partial_s \tilde x_{i}}{q_{i}}-
\frac{\partial_s \tilde x_{i}-\partial_s \tilde x_{i-1}}{q_{i-1}}
\right)
\\
&&
\approx
\left(\frac{x_{i+2}^j-x_{i+1}^j}{r_{i+1}^j q_{i+1}^j r_{i+2}^j}-
\frac{x_{i+1}^j-x_{i}^j}{r_{i+1}^j q_{i+1}^j r_{i+1}^j}\right)-
\left(\frac{x_{i+1}^j-x_{i}^j}{r_{i+1}^j q_{i}^j r_{i+1}^j}
-          
\frac{x_{i}^j-x_{i-1}^j}{r_{i+1}^j q_{i}^j r_{i}^j}\right)
\\
&&
-
\left(\frac{x_{i+1}^j-x_{i}^j}{r_{i}^j q_{i}^j r_{i+1}^j}-          
\frac{x_{i}^j-x_{i-1}^j}{r_{i}^j q_{i}^j r_{i}^j}\right)+
\left(\frac{x_{i}^j-x_{i-1}^j}{r_{i}^j q_{i-1}^j r_{i}^j}-              
\frac{x_{i-1}^j-x_{i-2}^j}{r_{i}^j q_{i-1}^j r_{i-1}^j}\right).
\end{eqnarray*}
The left hand side of (\ref{ode}) is approximated by means of backward time differences and by a central finite 
difference approximation of the advective term (up-wind technique can be also easily incorporated). We obtain
\[
q_i\frac{\hbox{d} x_i}{\hbox{d}t} +v_i (\tilde x_{i+1}-\tilde x_i)\approx
q_i^j\frac {x_{i}^j-x_i^{j-1}}{\tau} + v_i^j \frac {x_{i+1}^j-x_{i-1}^j}{2}\,.
\]
Taking into account information from the previous time step $j-1$ in the driving term $F \vec N,$
multiplying the third order terms approximation by $-\delta$ and putting all
unknowns $x_l^j$, $l=i-2,\dots,i+2$ with their coefficients to the left hand
side, we obtain a linear system of equations for updating the discrete position vector $x^j$ 
of the evolved curve $\Gamma^{j\tau}$:
\[
{\cal A}_{i}^j x_{i-2}^j+{\cal B}_{i}^j x_{i-1}^j+
{\cal C}_{i}^j x_{i}^j+ {\cal D}_{i}^j x_{i+1}^j+
{\cal E}_{i}^j x_{i+2}^j={\cal F}_i^j\
\]
for $i=1,...,n,$ subject to periodic boundary conditions
$x_{-1}^j=x_{n-1}^j,\ \ x_0^j=x_n^j,\ \ x_{n+1}^j=x_1^j,\ \ x_{n+2}^j=x_2^j$
where
\begin{eqnarray}
&&{\cal A}_i^j = \frac {\delta}{r_{i}^j q_{i-1}^j r_{i-1}^j},\ \ \  
{\cal E}_i^j = \frac {\delta}{r_{i+1}^j q_{i+1}^j r_{i+2}^j}\,,
\nonumber\\
&&{\cal B}_i^j = -\delta\left(\frac 1{r_i^j q_{i-1}^j r_{i-1}^j}
+\frac 1{(r_i^j)^2 q_{i-1}^j}+\frac 1{(r_i^j)^2 q_i^j}+
\frac 1{r_i^j q_i^j r_{i+1}^j}\right)
- \frac{\phi_i^j}{r_i^j} - \frac{v_{i}^j}{2}\,,
\nonumber\\
&&{\cal D}_i^j = -\delta\left(\frac 1{r_i^j q_{i}^j r_{i+1}^j}
+\frac 1{(r_{i+1}^j)^2 q_{i}^j}+\frac 1{(r_{i+1}^j)^2 q_{i+1}^j}
+\frac 1{r_{i+1}^j q_{i+1}^j r_{i+2}^j}\right)
-\frac{\phi_{i+1}^j}{r_{i+1}^j}
+\frac{v_{i}^j}{2}\,,
\nonumber\\
&&{\cal C}_i^j = \frac{q_i^j}{\tau}
-({\cal A}_i^j+{\cal B}_i^j+{\cal D}_i^j+{\cal E}_i^j)\,,\nonumber \\
&&{\cal F}_i^j = \frac {q_i^j}{\tau} x_i^{j-1}+q_i^j\
F(x_i^{j-1})\vec N\left(\frac{\nu_{i}^j+\nu_{i+1}^j}{2}\right)\,.
\nonumber
\end{eqnarray}
The previous system is penta-diagonal if $\delta>0$ and tri-diagonal
in the case $\delta=0$. 
In the latter case when $\delta=0$, the monotonicity assumption
on the function $b$ guarantees the strict diagonal dominance of the 
tri-diagonal system matrix.
In both cases, it is solved efficiently by means of the
Gauss-Seidel iterative method or by its well-known 
successive over relaxation version
SOR. We start the iterates 
from the previous time step vector $x^{j-1}$ and, in practice, when using
our asymptotically uniform tangential redistribution (AUTR), 
there are just few number of SOR 
iterations needed in order to achieve 
the solver accuracy goal. 
The iterative process is stopped when  a difference of subsequent iterates 
in the maximum norm is less than the prescribed tolerance, 
e.g. ${\rm TOL}=10^{-10}$.
Based on our practical experience, it should be also noted that the 
number of SOR iterations (corresponding to well-conditioning of an iteration matrix) 
is strongly lowered by such a proper choice of the tangential velocity. 
Consequently, it  speeds up computations
significantly making thus our numerical scheme fast and efficient.
Computational times are reported in next section.
In some examples of evolution of curves with high variation in the 
curvature or in checking the experimental order of convergence
(see e.g. Fig.~\ref{fig:spiralamcf} or Tables~\ref{tab:1} and \ref{tab:2}) 
we typically use small 
time steps $\tau\approx h^2$ 
or even $\tau\approx h^4$ (in some nontrivial cases of the fourth order flows).
Nevertheless, taking larger time steps $\tau$ leads to satisfactory 
numerical results as it can be seen from other experiments presented 
in this paper. We observed neither occurrences of accumulation of grid 
points nor spurious numerical oscillations. 

\section{Discussion on numerical experiments}

\begin{figure}
\begin{center}
\includegraphics[width=9cm]{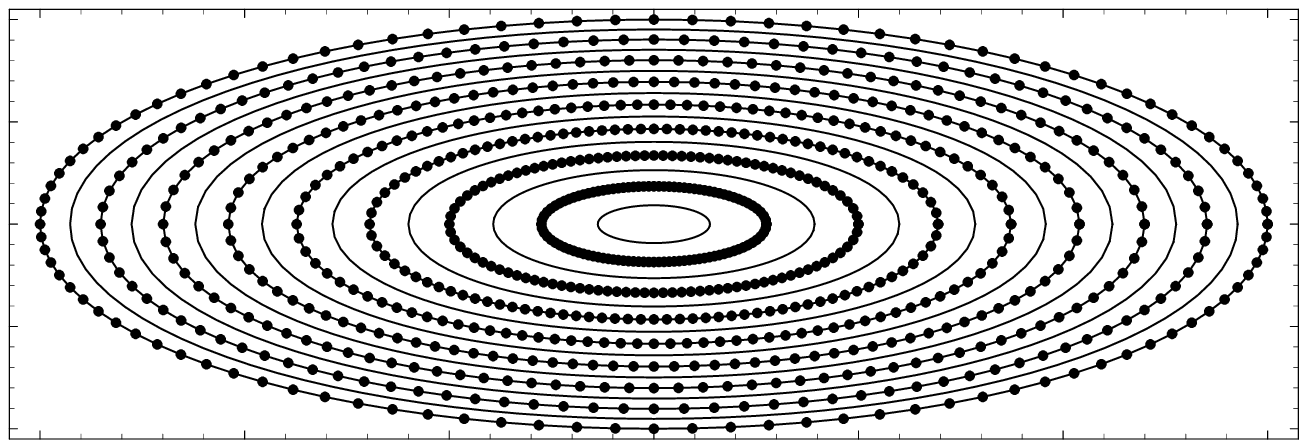}
\centerline{\scriptsize a)}
\vglue1mm
\includegraphics[width=9cm]{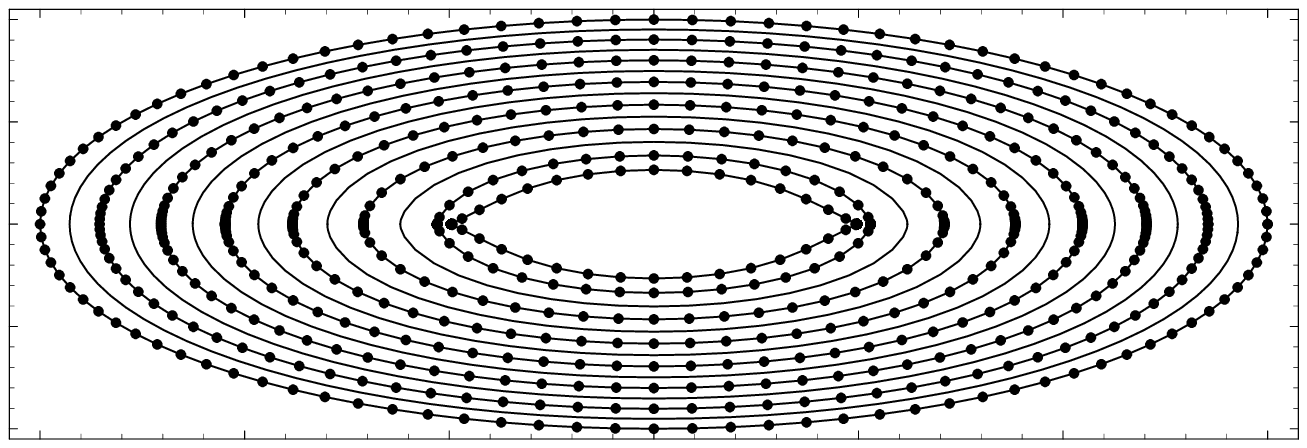}
\centerline{\scriptsize b)}
\end{center}
\caption{Affine invariant shrinking evolution of an initial ellipse 
($b=k^{1/3}, \delta=0, F= 0$); a) with asymptotically 
uniform redistribution ($\kappa_{1}= 3$); b) without redistribution. 
Numerical parameters $n=100, \tau = 0.001$. Time steps $j=0,  200, \cdots, 1400$ are plotted 
using  lines and grid points, while time steps $j= 100, 300, \cdots, $ up to $j=1500$ a); and $j=1100$ b) are plotted using lines.}
\label{fig:elipsa31beta13}
\end{figure}

\begin{figure}
\centering
\includegraphics[width=9cm]{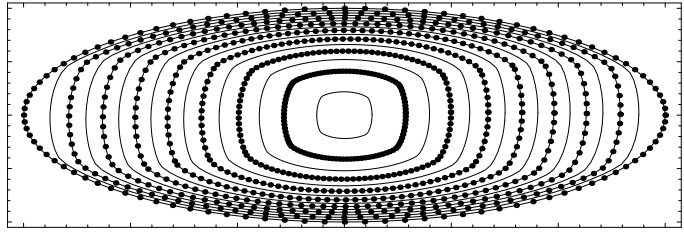}

\caption{Anisotropic shrinking evolution of an initial ellipse with AUTR;
$b=(1-0.9\cos(4\nu-\pi))\, k, \delta=0, F= 0$. Numerical and AUTR parameters: 
$n=100, \tau = 0.001$, $\kappa_{1}= 3$. Time steps $j=0,  200, \cdots, 1400$ are plotted 
using  lines and grid points, while time steps $j= 100, 300, \cdots,
1500$ are plotted using lines.}
\label{fig:anizotELIPSAredis}
\end{figure}

In the first numerical experiments we show a stabilizing effect of our scheme due to 
asymptotically uniform tangential redistribution (AUTR) 
for the case of selfsimilar affine invariant shrinking evolution of an initial ellipse with
half-axes ratio 3:1 (Fig.~\ref{fig:elipsa31beta13}a). When the grid
points are moving only in the normal direction, the numerical computation
collapses soon because of merging of grid points and spurious swallow-tails
creation in the left and right end of the ellipse (Fig.~\ref{fig:elipsa31beta13}b). 
In Fig.~\ref{fig:anizotELIPSAredis} we present an anisotropic curve shortening 
evolution of an ellipse computed again with help of the asymptotically uniform tangential redistribution.
In all numerical experiments parameters of computations are shown in figure captions.

\begin{figure}
\begin{center}
\includegraphics[width=3.5cm]{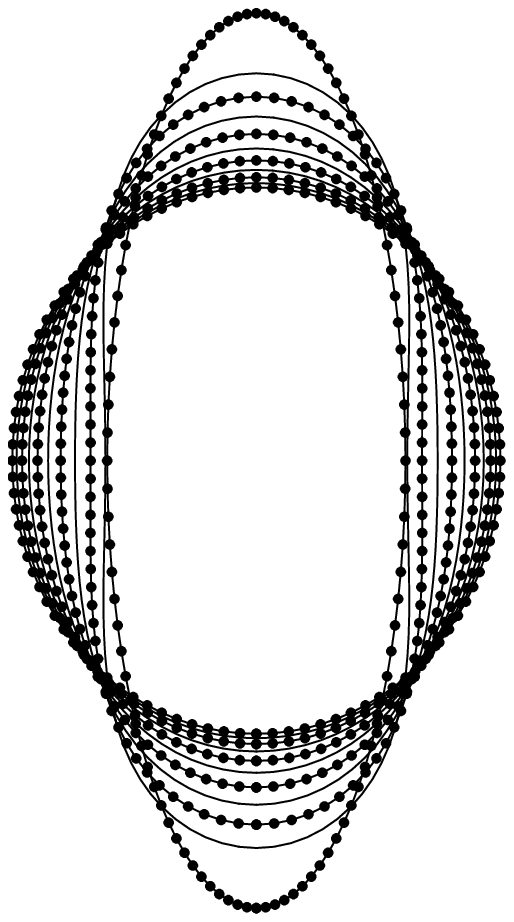}
\hglue1.5cm
\includegraphics[width=3.5cm]{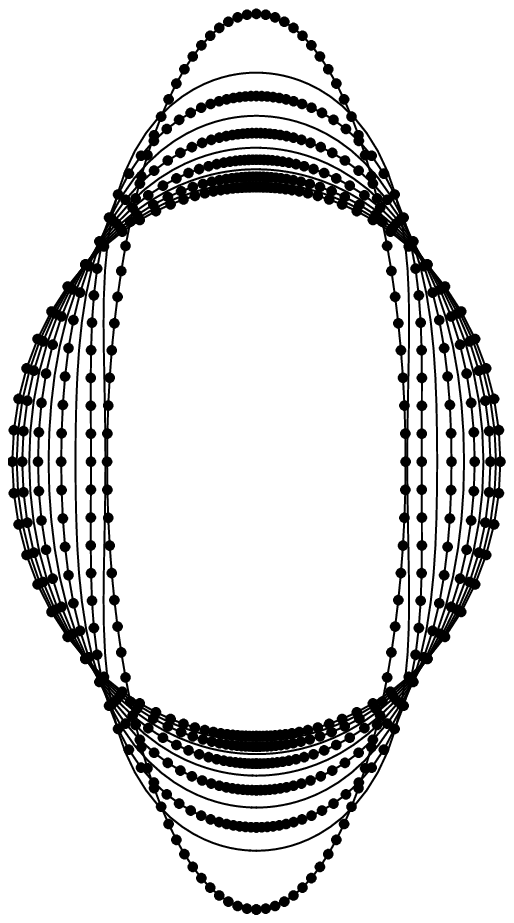}
\end{center}
\centerline{\scriptsize  a) \hglue 5truecm b)}
\caption{Surface diffusion flow ($\delta=1, b=0, F= 0$) of an initial ellipse 
with a) and without tangential redistribution b). 
Numerical and redistribution parameters: $n=100$, $\tau = 0.001$, 
$\kappa_{1}= 10$. Time steps $j=0, 400, \cdots, 2000$ are plotted
using  lines and grid points, while time steps $j= 200, 600, \cdots, 1800$ are plotted using lines.
}
\label{fig:elipsa31surfdif}
\end{figure}

\begin{table}
\caption{
An ellipse evolving by the surface diffusion using AUTR, same parameters as
in Fig.~\ref{fig:elipsa31surfdif}a. We report errors in area evolution, experimental
order of convergence  (EOC) in this quantity and computational time (CPU) 
for refined discretization parameters.}
\begin{center}
\footnotesize
\begin{tabular}{l|l|l|l|l|l} 

$n$& $\tau$&\# of steps & area error & EOC & CPU
(sec)\\ 
\hline
25  &  0.016   & 125 &  0.0774 &      & 0.02 \\
50  &  0.004   & 500 &  0.0202 & 1.93 & 0.43 \\
100 &  0.001   & 2000&  0.0052 & 1.89 & 2.94 \\
200 &  0.00025 & 8000&  0.0014 & 1.91 & 40.21\\
\end{tabular}
\end{center} 
\label{tab:1}
\end{table}

\begin{table}
\caption{
An ellipse evolving by the surface diffusion without tangential
redistribution, same parameters as
in Fig.~\ref{fig:elipsa31surfdif}b. 
Same quantities as above are reported.}
\begin{center}
\footnotesize
\begin{tabular}{l|l|l|l|l|l} 

$n$& $\tau$&\# of steps & area error & EOC & CPU
(sec)\\ 
\hline
25  &  0.016   & 125 &  0.2922 &      & 0.11 \\
50  &  0.004   & 500 &  0.1056 & 1.46 & 4.66 \\
100 &  0.001   & 2000&  0.0320 & 1.72 & 198 \\
200 &  0.00025 & 8000&  0.0155 & 1.04 & 1741\\
\end{tabular}
\end{center} 
\label{tab:2}
\end{table}

Again starting from an initial ellipse with half-axes ratio 1:3 
we numerically compute its evolution by the
surface diffusion. Fig.~\ref{fig:elipsa31surfdif}a  shows the result
with and Fig.~\ref{fig:elipsa31surfdif}b 
without asymptotically uniform
tangential redistribution (AUTR). In Tables~\ref{tab:1} and \ref{tab:2} we show 
how precisely the encompased area is preserved 
(which is one of the
analytical properties of surface diffusion flow with $\delta=1$ and $b=0$) 
during the evolution. 
The main
observation here is the fact that when using AUTR 
we obtain significantly lower error 
and the method with AUTR has approximately the second order of accuracy 
which is not the case for computations without AUTR.
In these experiments we used coupling $\tau\approx h^2$, cf.
\cite{D2,D3,DKS,MS_CVS,MS_ALG}, the 
exact area at any time moment is $A_e=3\pi$ and the area error is computed as
\begin{equation}
\label{aswee}
\left\|\epsilon^m_n \right\| =\biggl( \sum^{m}_{j=1}(A^{j} - A_{e})^2\ \tau \biggr)^\frac12 \,,
\end{equation}
where 
$A^j =  \frac{1}{2} \sum^{n}_{i=1} \det(x_i^j, x_i^j -x_{i-1}^j)$
gives the encompassed area of the polygonal computed curve at 
the $j-$th time step.

\begin{figure}%
\begin{center}
\includegraphics[width=5cm]{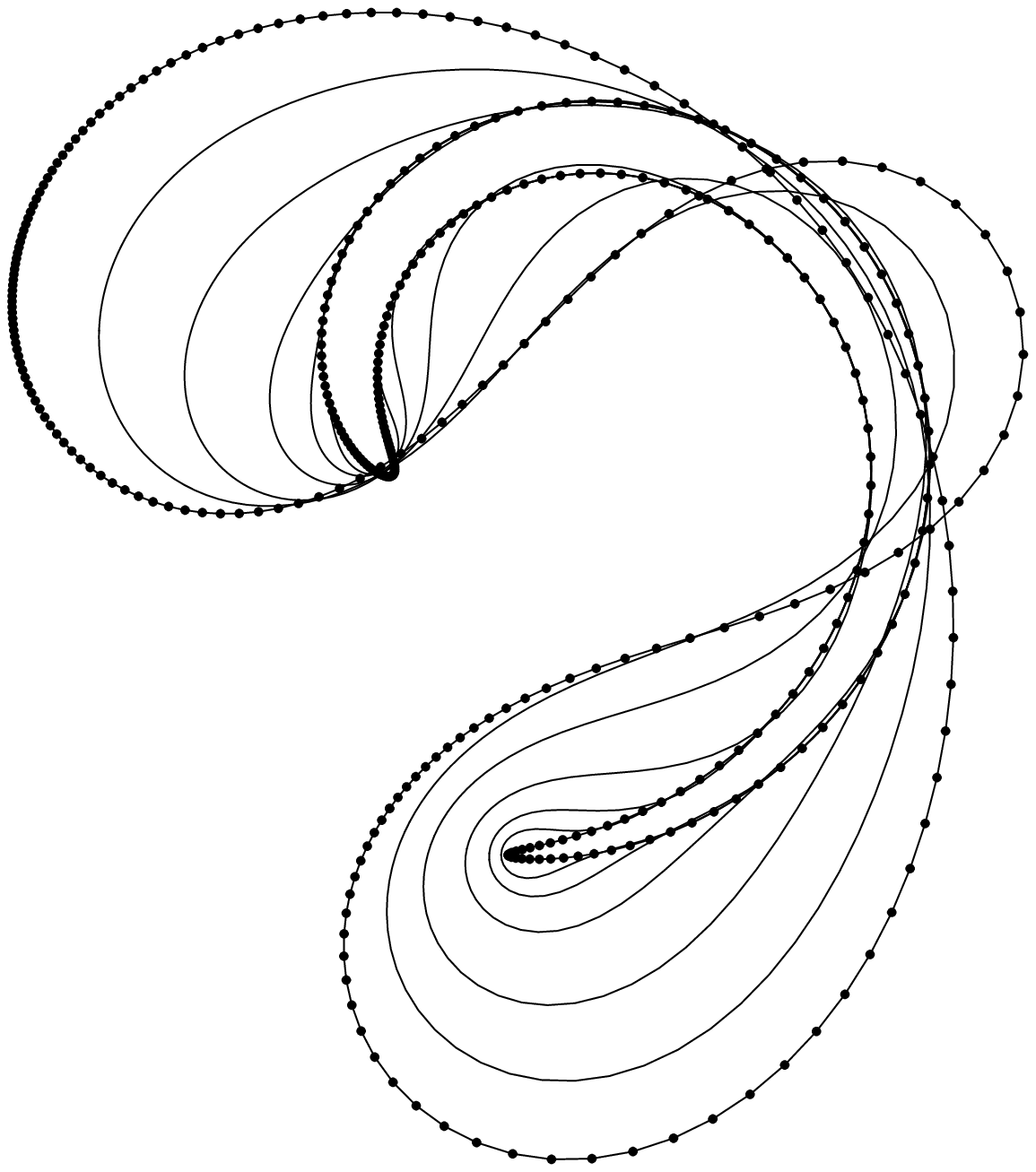}
\hglue3mm
\includegraphics[width=5cm]{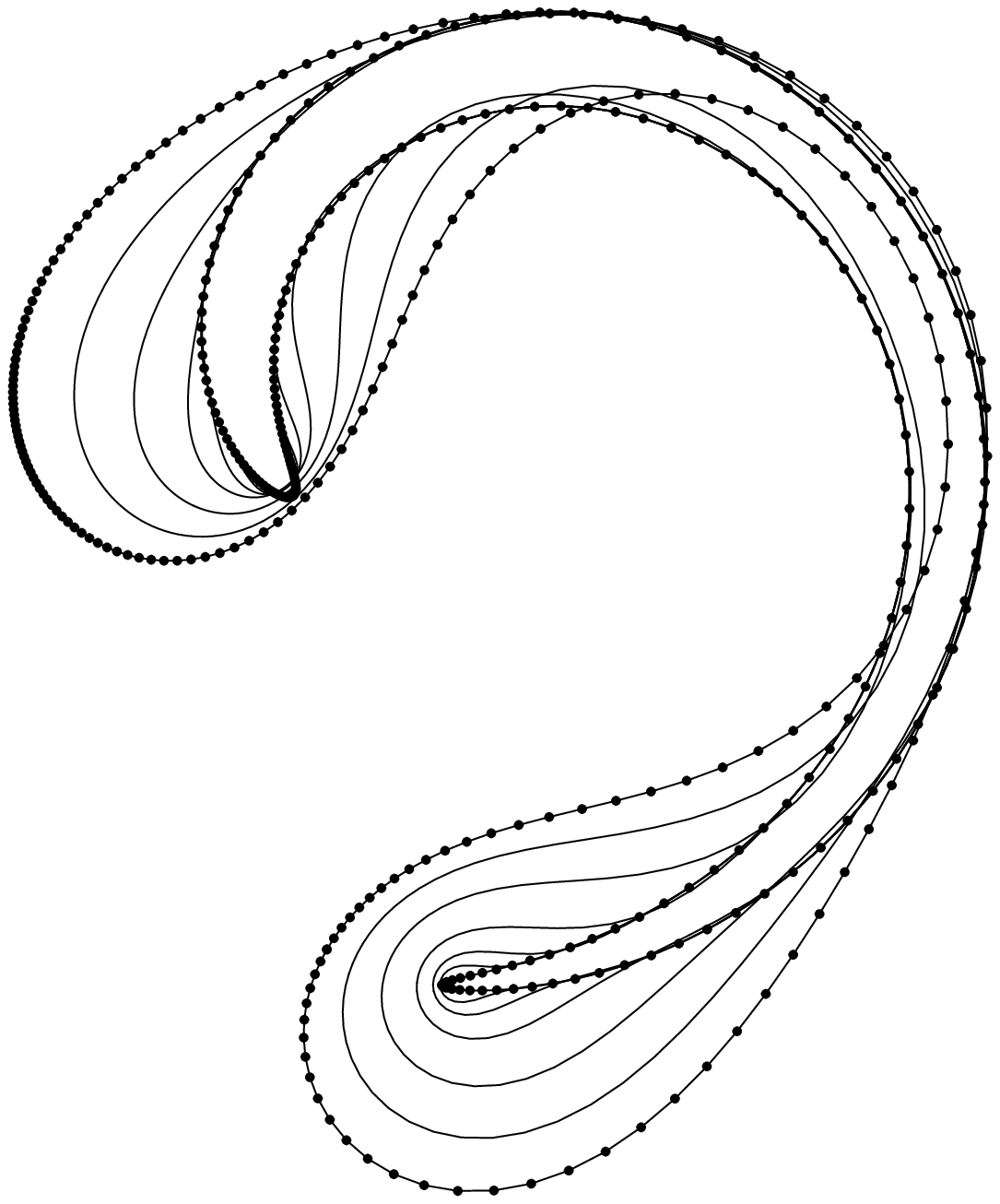}
\centerline{\scriptsize  a) \hglue6.5truecm b)}

\includegraphics[width=5cm]{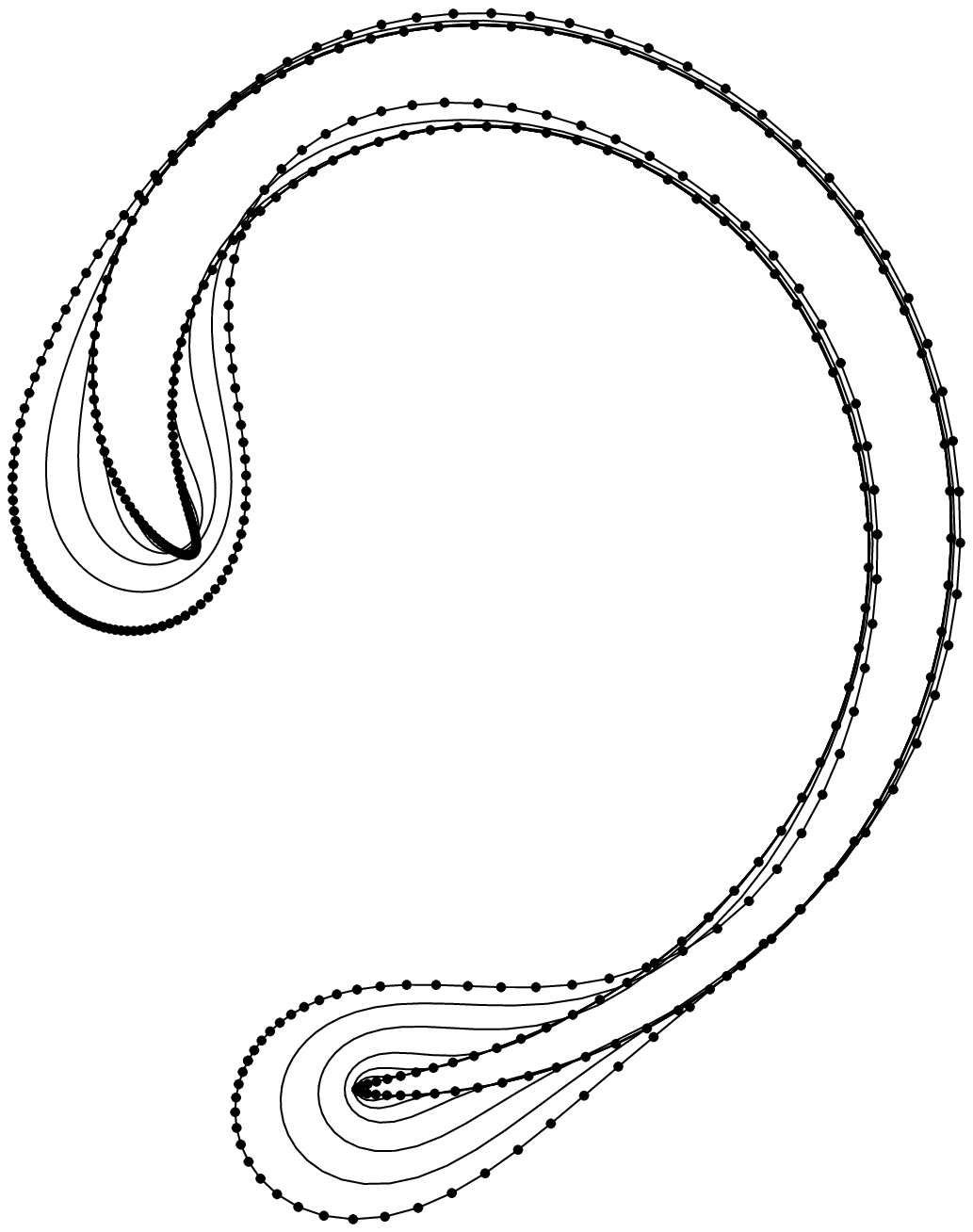}
\hglue3mm
\includegraphics[width=5cm]{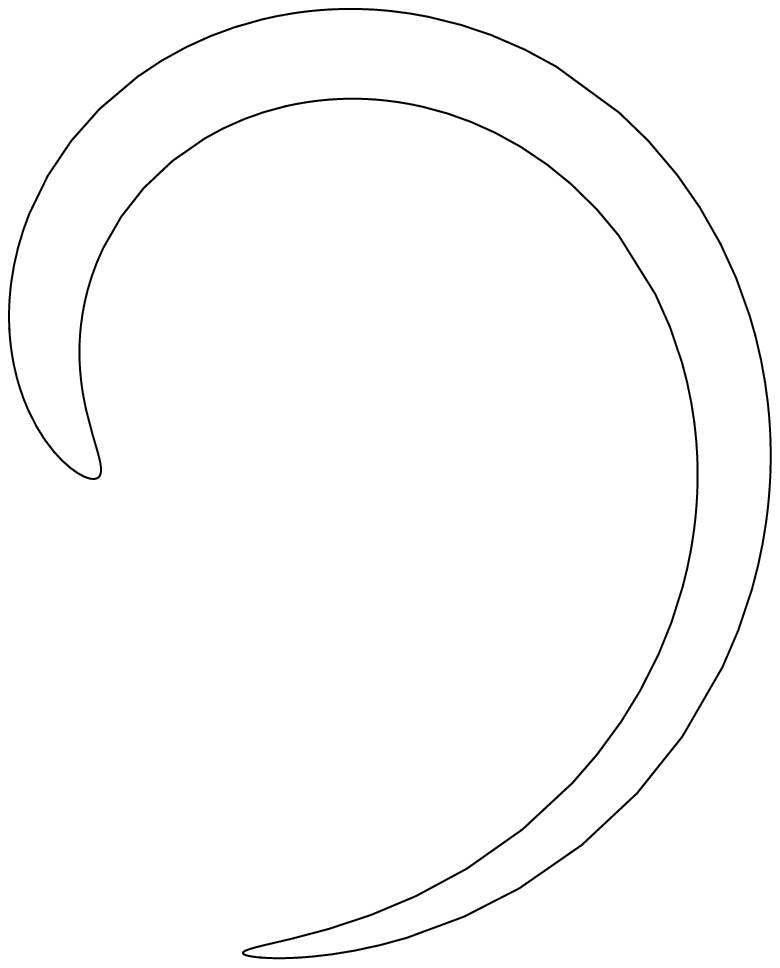}
\centerline{\scriptsize  c) \hglue6.5truecm d)}
\end{center}
\caption{Backward mean curvature flow regularized by the Willmore 
flow ($b(k,\nu)=-k-\frac 12 \delta k^3, F=0$) with different 
strengths $\delta=1$ a), $\delta=0.1$ b), $\delta=0.01$ c) 
of the same initial spiral d).
Numerical and AUTR parameters: $n=200, \tau = 10^{-12}$, $\kappa_{1}= 100$.
Time steps: $j=0, 10^5, 10^6, 10^7, 10^8,4\ 10^8,10^9$ are plotted. 
Initial and  the last step plotted with grid points.}
\label{fig:spiralaHAD}
\end{figure}

In Tables~\ref{tab:1} and \ref{tab:2} we also report computational times
achieved on a standard 2.2GHz laptop. They justify why the method with AUTR is refereed to as a {\it fast} method.
For standard curve resolutions, with 100 or 200 grid points, we get that one
time step takes 0.0014 respectively 0.0050 second. Such fast CPU times are 
obtained
due to tangential redistribution which not only stabilized the computations
but also improve cyclic penta-diagonal iteration matrix properties 
in the SOR iterative method (the relaxation parameter was set to $1.6$). Such CPU
times are obtained also in further experiments presented in this section. Of
course, one has to multiply them by the number of time steps which may lead to
overall long computations when computing long time behaviors of curve evolutions shown e.g. in Fig.~\ref{fig:spiralaHAD}.

\begin{figure}
\begin{center}
\includegraphics[width=6.45cm]{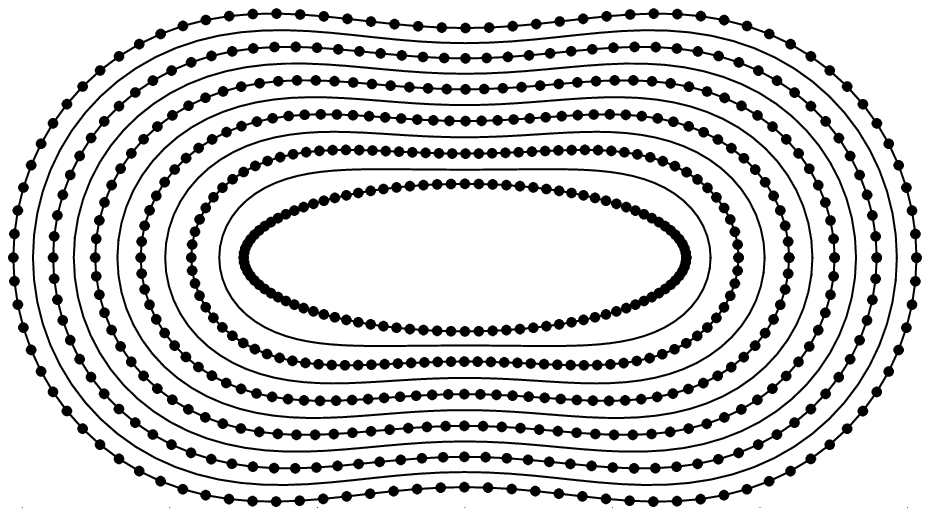}
\includegraphics[width=6.3cm]{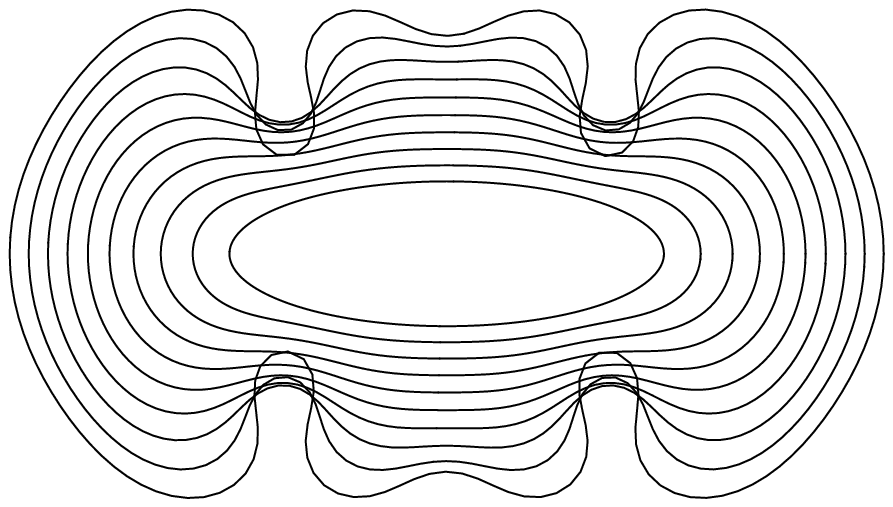}
\centerline{\scriptsize  a) \hglue6.5truecm b)}
\end{center}

\caption{Backward mean curvature flow with negative external force 
regularized by the Willmore surface energy 
($b=-k-\frac 12 \delta k^3, F=-1 $). a) a strong regularization effect 
with $\delta=1$;  b) a weak regularization $\delta=0.1$. 
Numerical and AUTR parameters:
$n=100$, $\tau = 0.001$, $\kappa_{1}= 10$. 
In subfigure a) the time steps $j=0, 400, \cdots, 2000$ are plotted using  
lines and grid points, while time steps $ j=200, 600, \cdots, 1800$ are 
plotted only with lines. In subfigure b) the 
time steps $j=0,  200, \cdots, 1800$ are plotted.}
\label{fig:elipsa31fm1mcfm1wfm1}
\end{figure}

\begin{figure}
\begin{center}
\includegraphics[width=6cm]{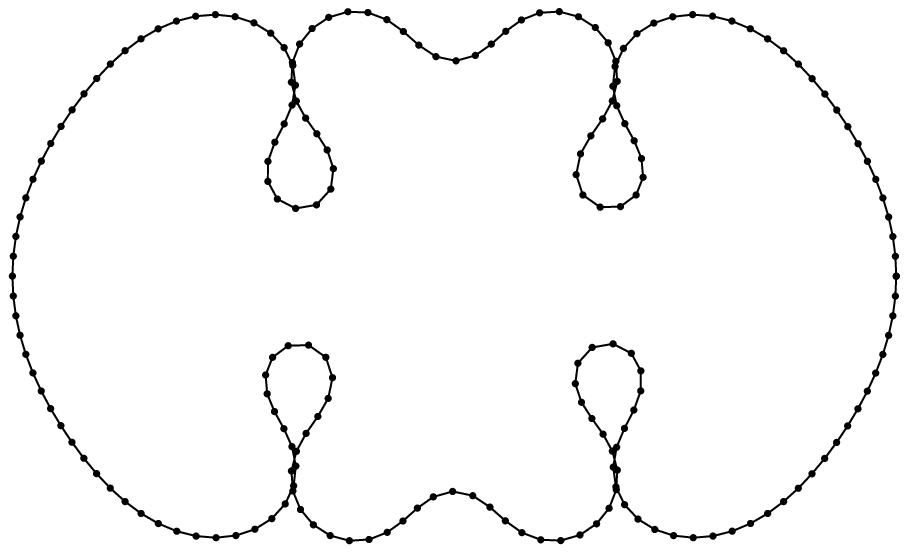}  
\includegraphics[width=6cm]{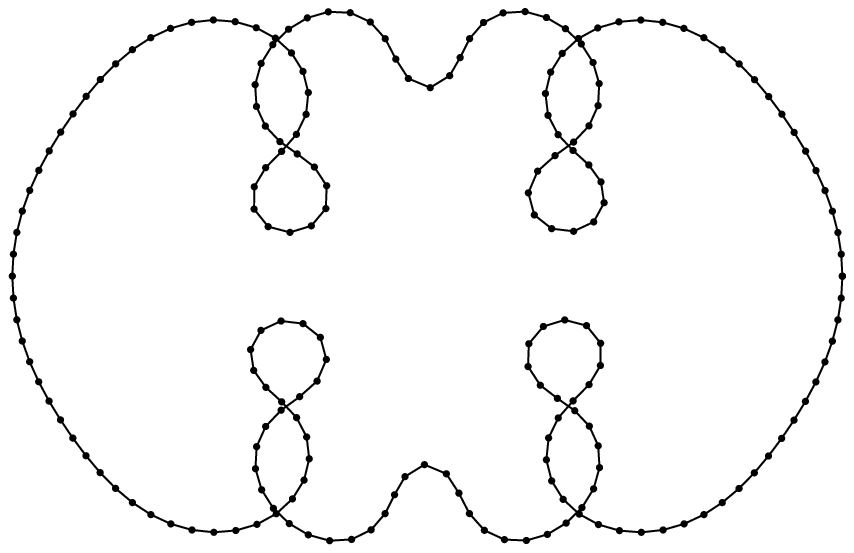} 
\end{center}
\begin{center}
\includegraphics[width=6cm]{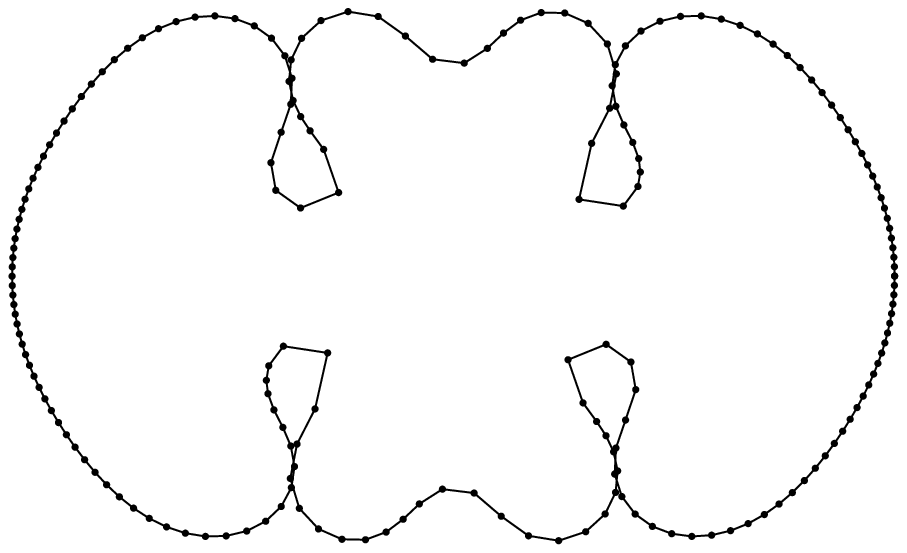}  
\includegraphics[width=6cm]{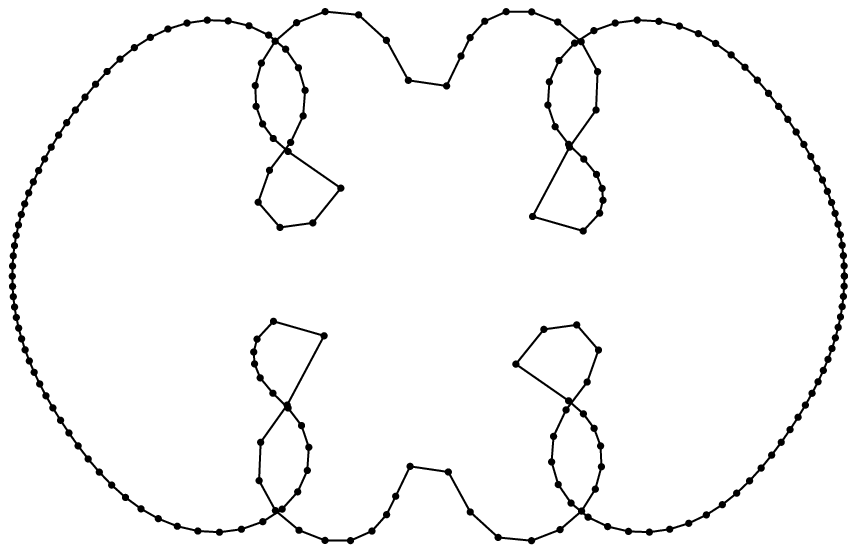}
\end{center}
\caption{Comparison of continuation of the evolution 
from Figure \ref{fig:elipsa31fm1mcfm1wfm1} b) 
computed with AUTR (upper row) and with redistribution preserving
relative local length (bottom row), respectively. In both cases we show
time steps  $j=2000$ (left) and $j=2200$ (right).}
\label{fig:elipsapokracovanie}
\end{figure}

\begin{figure}
\begin{center}
\includegraphics[width=8.5cm]{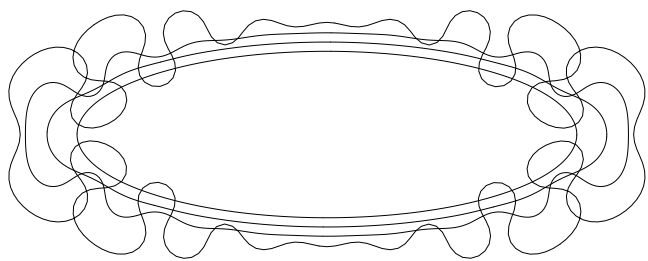}
\end{center}

\caption{Regularized backward curvature driven flow with $b=-k-\frac 12 
\delta k^3, F= -1$ and even weaker regularization with
$\delta=0.01$. Redistribution parameter: $\kappa_{1}= 10$ 
and numerical parameters: $n=400, \tau = 0.0001,$ are used. 
Time steps $j=0, 1000, 2000, 3000$ are plotted.}
\label{fig:elipsa31fm1mcfm1wf01b}
\end{figure}

The next set of experiments is focused on the backward mean curvature flow
with expanding constant force. The flow is regularized by various strengths of the 
Willmore flow. 
Starting with an initial ellipse with half-axes ratio 3:1 we can observe 
that in the case of a strong regularization the backward mean curvature flow 
is dominated by the elastic relaxation due to the Willmore surface energy 
(see Fig.~\ref{fig:elipsa31fm1mcfm1wfm1}a). On the other hand, if we decrease 
the fourth order regularization by taking smaller values of the parameter 
$\delta>0$, the nonconvex parts formed during evolution are attenuated
(see Fig.~\ref{fig:elipsa31fm1mcfm1wfm1}b) and even curve selfintersections 
may occur (see Figs.~\ref{fig:elipsapokracovanie} and 
\ref{fig:elipsa31fm1mcfm1wf01b}) as it can be expected in the backward 
in time diffusion process. 
 
In Fig.~\ref{fig:elipsapokracovanie} we also show an
important role of AUTR in such nontrivial experiments. The asymptotically
uniform tangential redistribution keeps very good curve resolution 
even in cases of several
subsequent curve selfintersections. On the other hand, if we consider for example the well-known redistribution preserving the relative local length \cite{Hou1,K2,MS2},
obtained by taking $\omega=0$ in
(\ref{uniformredis}) and (\ref{alpha-update}) then this method is not capable to handle this situations properly. 
If we consider  $\delta=0.01$ the backward diffusion effects are strongly dominating. 
In Fig.~\ref{fig:elipsa31fm1mcfm1wf01b} we can observe very fast 
shape "enrichment" soon after even a nonconvexification of the evolving curve.
It is worth to note that such experiments would be impossible without 
incorporating the tangential stabilization into the direct 
Lagrangian computational approach.

\begin{figure}
\begin{center}
		\includegraphics[width=5cm]{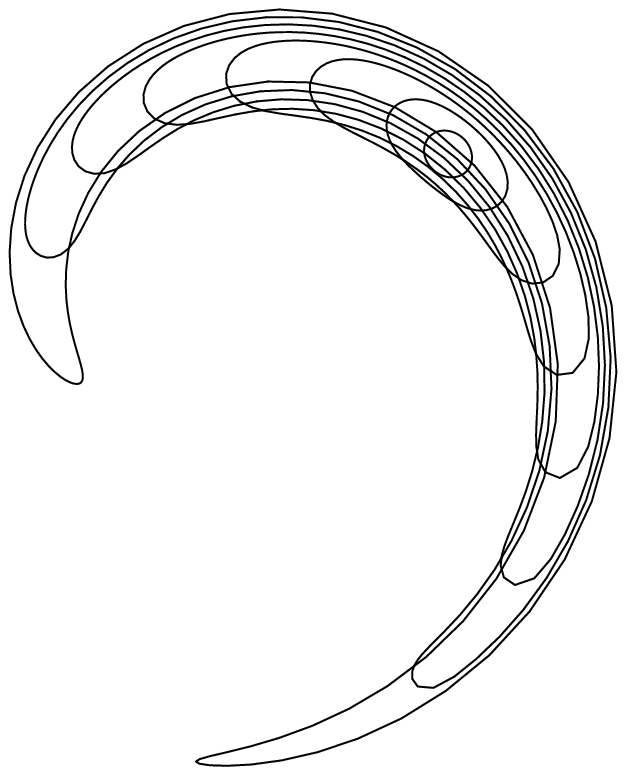}
\hglue 1truecm
		\includegraphics[width=5cm]{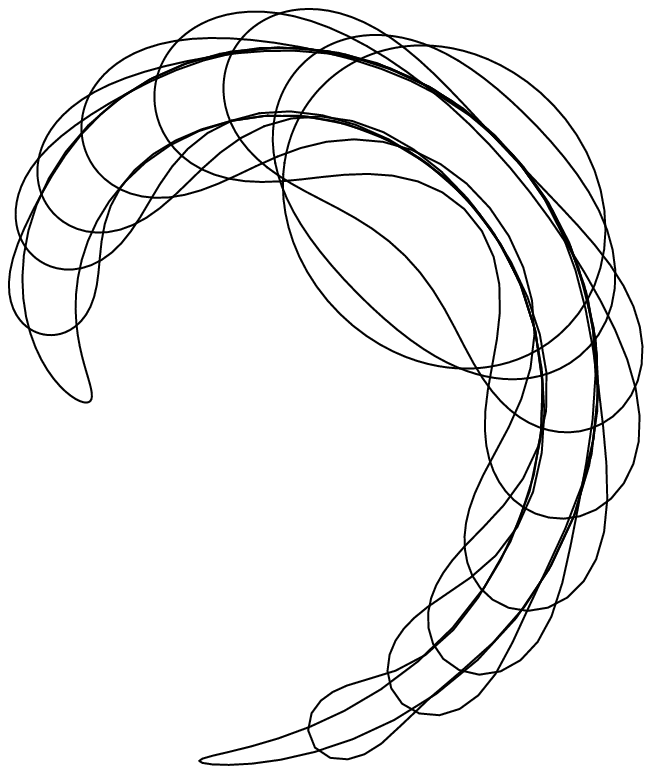}
\centerline{\scriptsize  a) \hglue6.5truecm b)}
\end{center}
\caption{a) mean curvature flow with $b(k,\nu)=k,\delta=0, F=0$; 
b) surface diffusion flow with $b=0, \delta=1, F=0$ of an initial spiral. 
Numerical and AUTR parameters: $n=100, \kappa_{1}= 10$ 
and $\tau = 10^{-6}$ a) and $\tau = 10^{-10}$ b). 
Time steps $j\in\{0, 1,2,3,4,5,6,6.5\} \times 10^4$ a) 
and $j\in\{0, 1,5,10,20,40,60,80,100\} \times 10^4$ b) are plotted.
}
\label{fig:spiralamcf}
\end{figure}

The last set of experiments is devoted to evolution of an initial spiral given by
\begin{eqnarray}
&&x_1(u)= a \cos b, \ \ \ x_2(u)= a \sin b\,,
\nonumber\\
&&a=0.5\ e^{-1-\frac 12 \sin(2\pi u)}-0.025\cos(2\pi u),\ \ \ 
b= 10\ \arctan(1+0.5 \sin(2\pi u))\,,
\nonumber
\end{eqnarray}
plotted in Fig.~\ref{fig:spiralaHAD}d.
Again the presence 
of the tangential redistribution  (AUTR) has a stabilizing effect on all numerical
computations. Without redistribution a parametric approach collapses soon. The
evolution by the mean curvature and surface diffusion is shown in 
Fig.~\ref{fig:spiralamcf}a and Fig.~\ref{fig:spiralamcf}b. 
The backward curve diffusion regularized  by a different strength 
of the Willmore flow is presented in Fig.~\ref{fig:spiralaHAD}a-c.

\section{Conclusions}
A new direct Lagrangian method stabilized by a suitable tangential 
redistribution has been presented for
the case of general plane curve evolution models. An evolved curve is driven 
in the normal direction
by a combination of the fourth order terms related to the intrinsic Laplacian of curvature,
second order terms related to the curvature, 
first order terms related to anisotropy and tangential redistribution 
and by a given external velocity field. We showed how a proper choice of a 
tangential velocity can stabilize and speed up computations. 
Nontrivial numerical experiments justified applicability and 
numerical stability of the approximation scheme in the 
anisotropic mean curvature flow, surface diffusion and the Willmore flow 
and even in the case of
backward curve diffusion slightly regularized by the fourth order terms.
Further study of the scheme from the numerical analysis point of view  as
well as derivation of analytical quantities to which the numerical results
can be compared in the above mentioned models will
be the objective of our future research.

\end{document}